%%%%%%%%%%%%%%%%%%%%%%%%%%%%%%%%%%%%%%%%%%%%%%%%%%%%%%%%%%%%%
%%                                                         %%
%% This file is the AMS-TeX file for the paper             %%
%%                                                         %%
%%   Real polynomials with all roots on the unit circle    %%
%%        and abelian varieties over finite fields         %%
%%                           by                            %%
%%         Stephen A. DiPippo and Everett W. Howe          %%
%%                                                         %%
%%                                                         %%
%% If you have problems with or questions about this file, %%
%% please contact me by email at                           %%
%%              however@alumni.caltech.edu                 %%
%%                                                         %%
%%%%%%%%%%%%%%%%%%%%%%%%%%%%%%%%%%%%%%%%%%%%%%%%%%%%%%%%%%%%%
\input amstex
\documentstyle{amsppt}

%% Define the AMS-TeX version 2.1 commands not found in version 2.0

\ifx\undefined\curraddr
  \def\curraddr#1\endcurraddr{\address {\it Current address\/}: #1\endaddress}
\fi

\ifx\undefined\rom
  \def\rom#1{{\rm #1}}
\fi

%% End of version 2.1 commands

%%%%%%%%%%%%%%%%%%%%%%%%%%%%%%%%%%%%%%%%%%%%%%%%%%%%%%%%%%%%%
%%                                                         %%
%%      These are the macros I use in the paper.           %%
%%                                                         %%
%%%%%%%%%%%%%%%%%%%%%%%%%%%%%%%%%%%%%%%%%%%%%%%%%%%%%%%%%%%%%

\define\BC{{\bold C}}
\define\BF{{\bold F}}
\define\BP{{\bold P}}
\define\BQ{{\bold Q}}
\define\BQpb{{\overline{{\bold Q}_p}}}
\define\BQp{{{\bold Q}_p}}
\define\BR{{\bold R}}
\define\BZ{{\bold Z}}

\define\bb{{\bold b}}
\define\bc{{\bold c}}
\define\be{{\bold e}}
\define\bm{{\bold m}}
\define\br{{\bold r}}
\define\bw{{\bold w}}
\define\bx{{\bold x}}
\define\By{{\bold y}}
\define\bz{{\bold z}}
\define\CI{{\Cal I}}
\define\CN{{\Cal N}}
\define\CO{{\Cal O}}

\let\hat=\widehat

\define\Chi{{\roman{X}}}

\define\Jac{\mathop{\roman{Jac}}\nolimits}
%\define\rad{\mathop{\roman{rad}}\limits}
%\define\res{\mathop{\roman{res}}\nolimits}
\define\volume{\mathop{\roman{volume}}\nolimits}
\define\covol{\mathop{\roman{covolume}}\nolimits}

\define\ra{\rightarrow}
%\define\la{\leftarrow}
%\define\ovl{\overline}
%\define\eps{\varepsilon}

% This next definition is used to give me my end-of-proof
% marks.  The \hfill is a harmless device used to put the 
% hollow square at the end of the line.  Please do not
% remove it --- end-of-proof marks *belong* at the end
% of the line (in my opinion).
\define\myqed{\hbox to 1em{}\nobreak\hfill$\square$}

% The numbers of the parts and sections of the paper...

\define\PartIntro{1}
\define\PartPoly{2}
\define\PartAV{3}

\define\SSregion{2.1}
\define\SSvolume{2.2}
\define\SSbounds{2.3}
\define\SSlipshitz{2.4}
\define\SSlowerbounds{2.5}

\define\SSweilpolys{3.1}
\define\SSisogenyclass{3.2}
\define\SSgrouporders{3.3}

% Theorem numbers...

\define\TTasymptotic{1.1}
\define\TTerrorterms{1.2}
\define\TTlowerbounds{1.3}
\define\TTgrouporders{1.4}

\define\LLhomeomorphism{2.1.1}
\define\EEregion2{2.1.2}
\define\LLboundary{2.1.3}

\define\PPvolume{2.2.1}
\define\PPvdMonde{2.2.2}

\define\PPlatticepoints{2.3.1}
\define\LLlatticepointsI{2.3.2}
\define\LLlipshitz{2.3.3}
\define\LLexactlattice{2.3.4}

\define\LLseverallipshitz{2.4.1}
\define\LLtwothirds{2.4.2}
\define\LLonehalf{2.4.3}

\define\LLniceregion{2.5.1}
\define\LLwedge{2.5.2}
\define\LLdiamond{2.5.3}

\define\PPordinary{3.1.1}
\define\LLnewton{3.1.2}
\define\PPnonordinary{3.1.3}

\define\PPAVlattice{3.2.1}
\define\LLAVdiamond{3.2.2}

\define\LLWeilpolys{3.3.1}
\define\EXgrouporders{3.3.2}

% ... and equation numbers.

\define\EQdefineChi{(1)}
\define\EQintegral{(2)}
\define\EQresidues{(3)}
\define\EQmessyinequality{(4)}
\define\EQbinsum{(5)}
\define\EQnumpoints{(6)}
\define\EQa1{(7)}
\define\EQai{(8)}

% The next macros are for citing references.

\define\refApostol{1}    % Apostol, Mathematical Analysis
\define\refGoss{2}       % Goss, Basic Structures of...
\define\refMarcus{3}     % Marcus, Number Fields
\define\refMehta{4}      % Mehta, Random Matrices (2nd ed)
\define\refOkada{5}      % Okada, On the generating functions...
\define\refRobbins{6}    % Robbins preprint
\define\refSelberg{7}    % Selberg, Bemerkinger...
\define\refTate{8}       % Tate, Classes d'isogenie...
\define\refWaterhouse{9} % Waterhouse's thesis

%%%%%%%%%%%%%%%%%%%%%%%%%%%%%%%%%%%%%%%%%%%%%%%%%%%%%%%%%%%%%
%%                                                         %%
%% If the version of the paper is not yet fixed, then      %%
%% it is today's date.                                     %%
%%                                                         %%
%%%%%%%%%%%%%%%%%%%%%%%%%%%%%%%%%%%%%%%%%%%%%%%%%%%%%%%%%%%%%

\define\version{20000421 --- see
J. Number Theory {\bf 73} {\rm(}1998\/{\rm)} 426--450\/{\rm;}\\
Corrigendum, J. Number Theory {\bf 83} {\rm(}2000\/{\rm)} 182.}

%\define\version{\number\year
%\ifcase\month\or 01\or 02\or 03\or 04\or 05\or 06\or 07\or 08\or 09%
%\or 10\or 11\or 12\fi
%\ifcase\day\or 01\or 02\or 03\or 04\or 05\or 06\or 07\or 08 \or 09\or 10%
%\or 11\or 12\or 13\or 14\or 15\or 16\or 17\or 18 \or 19\or 20%
%\or 21\or 22\or 23\or 24\or 25\or 26\or 27\or 28 \or 29\or 30%
%\or 31\fi \ --- under major revision}
 
% Tell TeX to let the user know what version he is TeXing. 
\message{Version \version}
 
% NOTE:  The \dedicatory macro below prints the version
% number on the TeX output.  

%%%%%%%%%%%%%%%%%%%%%%%%%%%%%%%%%%%%%%%%%%%%%%%%%%%%%%%%%%%%%
%%                                                         %%
%% Those are the last of my macros.  On to the topmatter...%%
%%                                                         %%
%%%%%%%%%%%%%%%%%%%%%%%%%%%%%%%%%%%%%%%%%%%%%%%%%%%%%%%%%%%%%

\topmatter
\title
Real polynomials with all roots on the unit circle
and abelian varieties over finite fields
\endtitle
\rightheadtext{REAL POLYNOMIALS AND ABELIAN VARIETIES}
\author
Stephen A. DiPippo and
Everett W. Howe
\endauthor
\address
Center for Communications Research,
29 Thanet Road,
Princeton, NJ 08540-3699 U.S.A.
\endaddress
\email
sdipippo\@idaccr.org
\endemail
\address
Center for Communications Research,
4320 Westerra Court,
San Diego, CA 92121-1967 U.S.A.
\endaddress
\email
however\@alumni.caltech.edu
\endemail

%%%%%%%%%%%%%%%%%%%%%%%%%%%%%%%%%%%%%%%%%%%%%%%%%%%%%%%%%%%%%
%%                                                         %%
%%  A kludge to get the version number on the preprint.    %%
%%                                                         %%
%%%%%%%%%%%%%%%%%%%%%%%%%%%%%%%%%%%%%%%%%%%%%%%%%%%%%%%%%%%%%
\dedicatory
Version \version
\enddedicatory 
% End of kludge.

%%%%%%%%%%%%%%%%%%%%%%%%%%%%%%%%%%%%%%%%%%%%%%%%%%%%%%%%%%%%%%%%%%%%%%%%%
%%                                                                     %%
%% 11          Number theory                                           %%
%%    11G         Arithmetic algebraic geometry (Diophantine geometry) %%
%%       11G10       Abelian varieties of dimension > 1                %%
%%       11G25       Varieties over finite and local fields            %%
%%                                                                     %%
%% 14          Algebraic geometry                                      %%
%%    14G         Arithmetic problems. Diophantine geometry            %%
%%       14G15       Finite ground fields                              %%
%%                                                                     %%
%% 33          Special functions                                       %%
%%    33B         Elementary classical functions                       %%
%%       33B15       Gamma, beta, and polygamma functions              %%
%%                                                                     %%
%%%%%%%%%%%%%%%%%%%%%%%%%%%%%%%%%%%%%%%%%%%%%%%%%%%%%%%%%%%%%%%%%%%%%%%%%


\abstract
In this paper we prove several theorems about abelian varieties over finite fields
by studying the set of monic real polynomials of degree $2n$
all of whose roots lie on the unit circle.  In particular,
we consider a set $V_n$ of vectors in $\BR^n$ that give the
coefficients of such polynomials. 
We calculate the volume of $V_n$
and we find a large easily-described subset of $V_n$.
Using these results, we find an asymptotic formula
--- with explicit error terms --- for the
number of isogeny classes of $n$-dimensional abelian varieties 
over $\BF_q$.  We also show that if $n>1$,
the set of group orders of $n$-dimensional abelian varieties
over $\BF_q$ contains every integer in an interval of length
roughly $q^{n-\frac{1}{2}}$ centered at $q^n+1$.
Our calculation of the volume of $V_n$ involves the evaluation
of the integral over the simplex 
$\bigl\{\,(x_1,\ldots,x_n) \bigm| 0\le x_1\le\cdots\le x_n\le 1\,\bigr\}$
of the determinant of the $n\times n$ matrix $\bigl[x_j^{e_i-1}\bigr]$, where the $e_i$ are
positive real numbers.
\endabstract

\endtopmatter

\head\PartIntro. Introduction
\endhead
In this paper we study the set of monic real polynomials of degree $2n$
with all roots on the unit circle, and we use our results on such polynomials
to prove several theorems about abelian varieties over finite fields. 
In particular, 
we find an asymptotic formula (for $n$ fixed and $q\ra \infty$)
for the number of isogeny classes of $n$-dimensional abelian varieties over $\BF_q$,
and we find a large interval of integers near $q^n+1$ that can be obtained as
the group orders of $n$-dimensional abelian varieties over $\BF_q$.

For every prime power $q$ and non-negative integer $n$  we let
$\CI(q,n)$ denote the set of isogeny classes of $n$-dimensional
abelian varieties over $\BF_q$; also,
we let $\CO(q,n)$ and $\CN(q,n)$ denote the ordinary and non-ordinary
isogeny classes in $\CI(q,n)$, respectively.
%Furthermore, we define rational numbers $v_n$ (for integers $n>0$) and an
%arithmetic function $r$ by setting 
%$$v_n = \frac{2^n}{n!}\prod_{j=1}^n \left(\frac{2j}{2j-1}\right)^{n+1-j}$$
%and $r(x) = \varphi(x)/x$, where $\varphi$ is Euler's $\varphi$-function.
Furthermore, for every positive integer $n$ we let 
$$v_n = \frac{2^n}{n!}\prod_{j=1}^n \left(\frac{2j}{2j-1}\right)^{n+1-j},$$
and we define an arithmetic function $r$ by setting
$r(x) = \varphi(x)/x$, where $\varphi$ is Euler's $\varphi$-function.

\proclaim{Theorem \TTasymptotic}
Let $n$ be a positive integer.  As $q\ra\infty$ over the prime powers,
we have $\#\CN(q,n) = O(q^{(n+2)(n-1)/4})$ and 
$\#\CI(q,n)\sim \#\CO(q,n) \sim v_n r(q) q^{n(n+1)/4}$.
\endproclaim

In fact, Theorem~\TTasymptotic\ is simply a corollary to a more
precise result that includes error terms.

\proclaim{Theorem \TTerrorterms}
For every positive integer $n$ and prime power $q$ we have
$$\left| \#\CO(q,n) - v_n r(q) q^{n(n+1)/4}\right|
 \le 6^{n^2} c_1^n c_2 \frac{n(n+1)}{(n-1)!} q^{(n+2)(n-1)/4}$$
and
$$ \#\CO(q,n-1) \le \#\CN(q,n) 
\le \left(v_n+6^{n^2} c_1^n c_3 \frac{n(n+1)}{(n-1)!}\right) q^{(n+2)(n-1)/4},$$
where 
$$\align
c_1 &= \sqrt{3}/6\approx 0.288675, \\
c_2 &= \exp(3/2) \cdot 2 (1+\sqrt{2}) \sqrt{3} (1+\sqrt{3}/162)^3/3 \approx 12.898608,\\
\vspace{-4\jot}
\intertext{and}
\vspace{-4\jot}
c_3 &= c_2 / (1+\sqrt{2}) \approx 5.342778.
\endalign$$
\endproclaim

If one is interested in obtaining lower bounds for the number 
of isogeny classes of $n$-dimensional abelian varieties over $\BF_q$
for specific values of $n$ and $q$, Theorem~\TTerrorterms\ is only
useful when $q$ is quite large compared to $n$; when $q$ is small
--- less than roughly~$6^{2n^2}$ --- 
one cannot even conclude from Theorem~\TTerrorterms\ that
$\#\CO(q,n)$ is nonzero.  To take care of this problem,
we prove a theorem that is not ideal asymptotically but that
does give nontrivial bounds when $q$ is small.

\proclaim{Theorem \TTlowerbounds}
For every positive integer $n$ and prime power $q$ we have
$$\#\CO(q,n) >  c_4 
(c_5 n)^{-2(\log 2)/(\log q)} \,\frac{2^n}{n!}  
\left(r(q)q^{n/2}-n\right) q^{n(n-1)/4}
,$$
where $c_4 = \exp(-3/2) \approx 0.223130$ and  $c_5 = 2 + \sqrt{2} \approx 3.414214$.
\endproclaim

Our other main topic concerns the group orders of abelian varieties over
finite fields.
If $A$ is an $n$-dimensional abelian variety over $\BF_q$, then 
Weil's ``Riemann Hypothesis'' 
shows that $\#A(\BF_q)$ is at least $(\sqrt{q}-1)^{2n}$ and 
at most $(\sqrt{q}+1)^{2n}$.
We prove that when $n>1$, every integer in a certain subinterval
of this allowable range actually does occur as the group order of an 
abelian variety.

\proclaim{Theorem~\TTgrouporders}
Let $q\ge 4$ be a power of a prime, 
let $$B_q = \frac{1}{2}\left(\frac{\sqrt{q}-2}{\sqrt{q}-1}\right) \text{\quad and\quad } 
C_q = \frac{\left\lfloor B_q\sqrt{q}\right\rfloor + 1/2}{\sqrt{q}},$$
and let $n>1$ be an integer.
If $m$ is an integer such that 
$\bigl|m - (q^n+1)\bigr| \le C_q q^{n-\frac{1}{2}}$
then there is an $n$-dimensional ordinary abelian variety $A$ over $\BF_q$
with $m = \#A(\BF_q)$.
\endproclaim

Note that the theorem would be false without the restriction that $n$
be greater than $1$, even if we were allowed to take $A$ to be 
non-ordinary; for example, if $q = p^a$ with $a>2$ then $q+1+p$ is
not the group order of an elliptic curve over $\BF_q$, as can be
seen from~\cite{\refWaterhouse, Theorem~4.1}.
Note also that when $q\gg n$ Weil's theorem 
restricts $\#A(\BF_q)$ to a range of roughly $4nq^{n-\frac{1}{2}}$ values, 
whereas about $q^{n-\frac{1}{2}}$ integers $m$ satisfy the inequality of Theorem~\TTgrouporders.

A result similar to Theorem~\TTgrouporders\ holds for $q <  4$  --- see
Exercise~\EXgrouporders.

The proofs of all of these theorems involve properties of the set $P_n$ of 
monic polynomials in $\BR[x]$ of degree $2n$ all of whose roots lie
on the unit circle and whose real roots occur with even multiplicity.
In Section~\SSvolume\ we find the volume of a region $V_n\subset\BR^n$ consisting
of vectors that give the coefficients of polynomials in $P_n$;
this volume computation involves the evaluation
of an integral that is reminiscent of the Selberg beta integral,
in that the integrand is the determinant of a Vandermonde-like matrix.
In Section~\SSbounds\ we give bounds (with error terms) on the sizes
of the intersections of $V_n$ with rectilinear lattices,
and in Section~\SSlowerbounds\ we find an easily-described subset
of $V_n$ that allows us to prove Theorems~\TTlowerbounds\ and~\TTgrouporders.

It is possible that our results on real polynomials with all roots on the
unit circle will have applications other than the ones we present here,
so we make almost no mention of abelian varieties in Part~\PartPoly\ of the
present paper.  We return to abelian varieties in Part~\PartAV,
where we prove the theorems in this introduction.

\subsubhead
Acknowledgments
\endsubsubhead
The authors thank: Joseph Wetherell, for asking questions that led to this research;
David Robbins, for suggesting Proposition~\PPvdMonde;
Hendrik Lenstra, for pointing out Lemma~\LLnewton;
Glenn Appleby, Noam Elkies, David Grabiner, and Richard Stanley,
for observing that certain special cases of the integral in 
Proposition~\PPvdMonde\ are special cases of the Selberg beta integral,
and for providing references; 
and Joshua Holden, for pointing out errors in the versions
of Theorem~\TTasymptotic\ and Proposition~\PPAVlattice\ that
appeared in the Journal of Number Theory.

\head\PartPoly. Real polynomials with all roots on the unit circle
\endhead

%%%%%%%%%%%%%%%%%%%%%%%%%%%%%%%%%%%%%%%%%%%%%%%%%%%%%%%%%%%%%%%%%%%%%%%%
\subhead\SSregion.  The parameterizing region $V_n$
\endsubhead

Suppose $g\in\BR[x]$ is a monic polynomial of degree $2n$ all of whose
complex roots lie on the unit circle, and suppose further that if $1$ 
or $-1$ is a root of $g$ then it occurs with even multiplicity.
Then the roots of $g$ come in complex-conjugate pairs, so $g$
factors over $\BR$
as the product of $n$ terms of the form $x^2 - r x + 1$, where $-2\le r\le 2$.
The symmetry of the factors of $g$ shows that there is a vector $\bb\in\BR^n$
such that $g=g_\bb$, where for every $\bb=(b_1,\ldots,b_n)$ we let $g_\bb$ denote the polynomial
$$g_\bb=\bigl(x^{2n} + 1\bigr) + b_1 \bigl(x^{2n-1}+x\bigr) + \cdots 
         + b_{n-1}\bigl(x^{n+1}+x^{n-1}\bigr) + b_n x^n.$$
We are led to the following definition:

\definition{Definition}
Let $V_n$ be the set of those $\bb\in\BR^n$ such that
all of the complex roots of $g_\bb$ lie on the unit circle.
\enddefinition

Our main goal in this section is to find an explicit homeomorphism $\Phi$ from a
simplex to $V_n$; this homeomorphism will be essential for our calculation
of the volume of $V_n$ in Section~\SSvolume.
The simplex that will be most convenient for us is the subset $I_n$ of $\BR^n$
defined by
$I_n=\bigl\{\,(x_1,\ldots,x_n) \bigm| -2 \le x_1 \le x_2 \le \cdots \le x_n \le 2\,\bigr\}$,
and the homeomorphism $\Phi$ will be obtained as the composition of two maps $\Psi$ and
$\Chi$ from $\BR^n$ to $\BR^n$ that we now define.

For every $\bc=(c_1,\ldots,c_n)\in \BR^n$ we let $h_\bc$ denote the polynomial
$h_\bc= x^n + c_1x^{n-1} + \cdots + c_n$.
If $\br = (r_1,\ldots,r_n)$ is an element of $\BR^n$, we let $\Psi(\br)$
be the vector $\bc\in\BR^n$ such that
$(x-r_1)\cdots(x-r_n) = h_\bc$
is an equality of polynomials; in other words, the $i$th component of $\bc$
is $(-1)^i$ times the $i$th symmetric polynomial in the~$r_j$.
Next, given a vector $\bc\in\BR^n$, we define $\Chi(\bc)$ 
to be the unique $\bb\in\BR^n$ such that 
$$x^n h_\bc(x+1/x) = g_\bb\eqno{\EQdefineChi}$$
is an equality of polynomials.  
Finally, we let $\Phi = \Chi\circ\Psi$.  

\proclaim{Lemma~\LLhomeomorphism}
The map $\Phi$ induces a homeomorphism from $I_n$ to $V_n$.
\endproclaim

\demo{Proof}
Suppose $\bb$ is an element of $V_n$. The symmetry of the polynomial $g_\bb$
shows that we may write 
$g_\bb = x^n h_\bc(x+1/x)$ for some unique $\bc\in\BR^n$, and it is clear from this equality
that the roots of $h_\bc$ are exactly
the traces (from $\BC$ to $\BR$) of the roots of~$g_\bb$.  Thus the roots of 
$h_\bc$ all lie between $-2$ and $2$, and it follows that $\bc=\Psi(\br)$ for
a unique vector $\br$ in $I_n$ --- namely, the vector that consists of the roots
of $h_\bc$ listed in non-decreasing order.  Conversely, if $\br\in I_n$ it is clear
that $\Chi\bigl(\Psi(\br)\bigr) \in V_n$.  Thus, $\Phi$ gives a bijection between $I_n$ and $V_n$.
The lemma then follows from the fact that
a continuous bijection from a compact topological space to a Hausdorff space
is a homeomorphism.
\myqed
\enddemo

\example{Example \EEregion2}
The set $V_1$ is just the closed interval $[-2,2]$,
and $V_2$ is also easily-described:
A simple calculation shows that for $n=2$
the map $\Phi\colon \BR^2\ra\BR^2$ sends $(r_1,r_2)$ to $(b_1,b_2)=(-r_1-r_2,r_1r_2+2)$,
and by noting where $\Phi$ sends the boundary of $I_2$, one checks that $V_2$ is the set
of all $(b_1,b_2)$ such that $b_2 \ge 2|b_1|-2$ and $b_2\le b_1^2/4 + 2$.
We note for future reference that the set $V_2$ contains the square $[-1,1]\times[0,2]$.
\endexample

Lemma~\LLhomeomorphism\ allows us to characterize the elements of the boundary
of $V_n$ in terms of the roots of their associated polynomials.

\proclaim{Lemma~\LLboundary}
The boundary of $V_n$ consists of those $\bb\in V_n$ 
such that $g_\bb$ has multiple roots.
\endproclaim

\demo{Proof}
Equation~\EQdefineChi\ tells us that for every $\br\in I_n$ 
the multi-set of roots of $g_{\Phi(\br)}$
is equal to the union over $i$ of the
multi-sets of roots of $x+1/x = r_i$, so $g_{\Phi(\br)}$ will have multiple roots
exactly when either two of the $r_i$ are equal to one another or one of the $r_i$ is $\pm 2$.
But the latter conditions are exactly the conditions for an element $\br$ of $I_n$ to be
in $\partial I_n$, so the set of $\bb\in V_n$ such that $g_\bb$ has
multiple roots is equal to $\Phi(\partial I_n)$, and this last set is $\partial V_n$
because $\Phi$ is a homeomorphism of manifolds-with-boundary.
\myqed
\enddemo

%%%%%%%%%%%%%%%%%%%%%%%%%%%%%%%%%%%%%%%%%%%%%%%%%%%%%%%%%%%%%%%%%%%%%%%%
\subhead\SSvolume.  A volume calculation
\endsubhead

In this section we will calculate the volume of the region $V_n$.

\proclaim{Proposition \PPvolume}
The volume $v_n$ of the region $V_n$ is given by
$$v_n = \frac{2^{n(n+1)}}{n!}\prod_{1\le i< j\le n}\frac{j-i}{j+i}
 = \frac{2^n}{n!}\prod_{j=1}^n \left(\frac{2j}{2j-1}\right)^{n+1-j}.$$
\endproclaim

\demo{Proof}
We maintain the notation used in Section~\SSregion.   Recall that if
$\bc = (c_1,\ldots,c_n)$ is an element of $\BR^n$ and if we let $\bb = \Chi(\bc)$,
then we have $x^n h_\bc(x+1/x) = g_\bb$, 
where $h_\bc$ and $g_\bb$ are the polynomials defined above.
Equating the coefficients in these polynomials, we find that
for every $i$ we have $$b_i = c_i + \hbox{(terms involving only $c_1,\ldots,c_{i-1}$)}.$$
Thus the Jacobian matrix $\Jac \Chi$ of $\Chi$ is triangular with $1$'s on the diagonal,
so $\det \Jac\Chi = 1$ and the volume of $V_n$ is equal to the volume
of $\Psi(I_n)$.

Let $\Omega\colon\BR^n\ra\BR^n$ denote the power-sum map
that sends $\br=(r_1,\ldots,r_n)$ to $(p_1,\ldots,p_n)$,
where $p_i = r_1^i + \cdots + r_n^i$.
Every $p_i$ is a symmetric polynomial in the variables $r_j$ and is
therefore expressible as a polynomial function of the elementary symmetric polynomials
in the $r_j$.  This implies that we have $\Omega = \Upsilon\circ\Psi$
for some polynomial map $\Upsilon\colon\BR^n\ra\BR^n$, because the coordinates of 
$\Psi$ are the elementary symmetric polynomials in the $r_i$ (up to sign).
Explicitly, if we let $\bc = \Psi(\br)$
then  Newton's formulas  relating the
power-sums to the symmetric polynomials state that
$$i c_i + \sum_{j=0}^{i-1} c_j p_{i-j} = 0$$
for all $i$; here we have set $c_0 = 1$.  From these equalities it follows that
for every $i$ we have
$$p_i = -i c_i + \hbox{(terms involving only $c_1,\ldots,c_{i-1}$)}.$$
Thus $\Jac\Upsilon$ is a triangular matrix whose $i$th diagonal entry
is $-i$, so $\Jac\Upsilon$ has determinant $(-1)^n n!$, and we see
that the volume of $\Psi(I_n)$ is equal to $1/n!$ times the
volume of $\Omega(I_n)$.

Now, since $\Omega$ is injective on $I_n$ the volume of $\Omega(I_n)$ is 
equal to the integral of $\bigl|\det\Jac\Omega\bigr|$ over $I_n$,
and $\det\Jac\Omega = \det\left[i x_j^{i-1}\right] = n!\, \det\left[x_j^{i-1}\right]$.
This Vandermonde determinant is always non-negative on $I_n$ so we may drop
the absolute value signs in the integrand.
Also, since the determinant depends only on the differences
of the $x_i$, we may shift the region of integration from $I_n$
to 
$$\bigl\{\,(x_1,\ldots,x_n) \bigm| 0 \le x_1 \le x_2 \le \cdots \le x_n \le 4\,\bigr\}$$
without changing the value of the integral.
If we then scale the variables by setting $x_i = 4 y_i$, we change
the region of integration to 
$$J_n = \bigl\{\,(y_1,\ldots, y_n)\in\BR^n \bigm| 0\le y_1\le\cdots\le y_n\le 1\,\bigr\}$$
at the expense of adding a factor
of $4$ for each of the $n\choose 2$ factors of the Vandermonde determinant
and each of the $n$ differentials $dx_i$.  We find that
$$\int_{I_n} \det\Jac\Omega\ dx_1\cdots dx_n
    = n!\, 4^{{n\choose2}+n}\int_{J_n} \det\left[y_j^{i-1}\right] \ dy_1\cdots dy_n,$$
and by using  Proposition~\PPvdMonde\ below we see that we have
$$\int_{I_n} \det\Jac\Omega\ dx_1\cdots dx_n
    = 2^{n(n+1)}\prod_{1\le i<j\le n}\frac{j-i}{j+i}.$$
Combining this equality with the conclusions of the preceding two
paragraphs, we find the first expression for~$v_n$ given in the
proposition.  A simple induction gives us the second expression as well.
\myqed
\enddemo

\proclaim{Proposition~\PPvdMonde\ (Robbins)}
Let $e_1$, \dots,~$e_n$ be positive real numbers and let $J_n$ be the
simplex define above.  Then 
$$\int_{J_n} \det\left[x_j^{e_i-1}\right] \ dx_1\cdots dx_n = 
\frac{1}{e_1\cdots e_n}\prod_{1\le i < j\le n} \frac{e_j-e_i}{e_j+e_i}.
\eqno{\EQintegral}$$
\endproclaim

\remark{Remark}
Robbins \cite{\refRobbins} provides an elegant proof of Proposition~\PPvdMonde\ based
on a result of Okada \cite{\refOkada, Theorem~3}.  
We present an alternate proof
that avoids the use of Okada's theorem.
\endremark
\remark{Remark}
We will only have call to apply 
Proposition~\PPvdMonde\ in the special case where $e_i = i$, in which
case the integral in the proposition is a special case of the
Selberg beta integral (see \cite{\refMehta, Chapter~17} or \cite{\refSelberg}).
\endremark

\demo{Proof of Proposition~{\rm \PPvdMonde}}
Let $D_n(e_1,\cdots,e_n)$ denote the left-hand side of~\EQintegral\ and
let $E_n(e_1,\cdots,e_n)$ denote the right-hand side.  
We will prove the following three statements, from which the
proposition follows:
\roster
\item"1." For every $e_1>0$ we have $D_1(e_1) = E_1(e_1)$.
\item"2." For $n>1$ we have
$$D_n(e_1,\ldots,e_n)  = \frac{1}{e_1+\cdots+e_n}
\sum_{k=1}^n (-1)^{n-k} D_{n-1}(e_1,\ldots,\hat{e_k},\ldots, e_n),$$
where a hat over a variable means that it is to be omitted.
\item"3." For $n>1$ we have 
$$E_n(e_1,\ldots,e_n)  = \frac{1}{e_1+\cdots+e_n}
\sum_{k=1}^n (-1)^{n-k} E_{n-1}(e_1,\ldots,\hat{e_k},\ldots, e_n).$$
\endroster

{\it Proof of statement\/} 1.
It is easy to check that $D_1(e_1)$ and $E_1(e_1)$ both equal $1/e_1$.

{\it Proof of statement\/} 2.
We will evaluate $D_n(e_1,\ldots,e_n)$ by integrating with respect to $x_n$
on the outside and rescaling the remaining variables by setting
$y_i = x_i/x_n$ for $i=1$, \dots,~$n-1$.  We find
$$D_n(e_1,\ldots,e_n) = \int\limits_{0\le x_n\le 1} \quad \int\limits_{0\le y_1\le\cdots\le y_{n-1}\le 1}
(\det M) x_n^{n-1} \ dy_1\cdots dy_{n-1}\ dx_n$$
where $M$ is the matrix 
$$M = 
\left[\matrix
x_n^{e_1-1}y_1^{e_1-1}  & \cdots & x_n^{e_1-1}y_{n-1}^{e_1-1} & x_n^{e_1-1}\\
    \vdots              & \ddots & \vdots                     & \vdots\\
x_n^{e_n-1}y_1^{e_n-1}  & \cdots & x_n^{e_n-1}y_{n-1}^{e_n-1} & x_n^{e_n-1}
\endmatrix\right].
$$
Every entry in the $i$th row of $M$ has a factor of $x_n^{e_i-1}$, so 
the determinant of $M$ is equal to 
$x^{(e_1-1) + \cdots + (e_n-1)}$ times the determinant of the matrix
$$N = 
\left[\matrix
y_1^{e_1-1}  & \cdots & y_{n-1}^{e_1-1} & 1\\
    \vdots   & \ddots & \vdots          & \vdots\\
y_1^{e_n-1}  & \cdots & y_{n-1}^{e_n-1} & 1
\endmatrix\right].
$$
Thus
$$\align
D_n(e_1,\ldots,e_n) & = \int\limits_{0\le x_n\le1} x_n^{e_1+\cdots+e_n-1}\ dx_n \quad
\int\limits_{0\le y_1\le\cdots\le y_{n-1}\le 1}
(\det N) \ dy_1\cdots dy_{n-1}\\
&= \frac{1}{e_1+\cdots+e_n} \int_{J_{n-1}} (\det N) \ dy_1\cdots dy_{n-1}.
\endalign
$$
Expanding $\det N$ by minors on the $n$th column gives us statement 2.

{\it Proof of statement\/} 3. 
{}From the definition of $E_n$ we see that 
$$\align
(-1)^{n-k} \frac{E_{n-1}(e_1,\ldots,\hat{e_k},\ldots, e_n)}{E_n(e_1,\ldots,e_n)}
& = (-1)^{n-k} e_k \Bigg(\prod_{i=1}^{k-1} \frac{e_k + e_i}{e_k-e_i}\Bigg)
                   \Bigg(\prod_{j=k+1}^{n} \frac{e_j + e_k}{e_j-e_k}\Bigg)\\
& = e_k \prod_{i\neq k} \frac{e_k + e_i}{e_k-e_i},
\endalign
$$
so to prove statement~3 it will be enough to show that
$$e_1+\cdots+e_n = \sum_{k=1}^n e_k \prod_{i\neq k} \frac{e_k + e_i}{e_k-e_i}.\eqno{\EQresidues}$$
Consider the meromorphic differential $\omega = \prod_{i=1}^n \bigl( (z+e_i)/(z-e_i)\bigr)\ dz$
on $\BP^1$, whose only poles are at the $e_k$ and at infinity.  It is easy to calculate
that the residue of $\omega$ at $z = e_k$ is
$2e_k \prod_{i\neq k} \bigl( (e_k + e_i)/(e_k-e_i)\bigr)$
and that the residue of $\omega$ at $z=\infty$ is
$-2(e_1+\cdots+e_n)$.
Equation~\EQresidues\ then follows from the fact that
the sum of the residues of a meromorphic differential on $\BP^1$ is zero.
\myqed
\enddemo

%%%%%%%%%%%%%%%%%%%%%%%%%%%%%%%%%%%%%%%%%%%%%%%%%%%%%%%%%%%%%%%%%%%%%%%%
\subhead\SSbounds.  Lattice points in $V_n$
\endsubhead

In the course of proving our asymptotic formulas for the number
of isogeny classes of $n$-dimensional abelian varieties over a
finite field, we will require estimates for the sizes of the
intersections of various lattices in $\BR^n$ with the set $V_n$ defined in 
Section~\SSregion.  In this section we will provide such estimates,
but only for the type of lattices we will encounter in our later
work: {\it rectilinear\/} lattices, by which we mean lattices that
have a rectilinear fundamental domain with edges parallel to the
coordinate axes.
The {\it covolume\/} of a lattice $\Lambda$ is the volume of 
a fundamental domain $R$ of $\Lambda$, and if $R$ is rectilinear
then the {\it mesh\/} of $\Lambda$ is the length of the longest 
edge of $R$.

\proclaim{Proposition~\PPlatticepoints}
Let $n>0$ be an integer and let $\Lambda\subset\BR^n$ be
a rectilinear lattice with mesh $d$ at most $1$.
Then we have
$$\left| \#(\Lambda\cap V_n) - \frac{\volume V_n}{\covol \Lambda} \right|
 \le 6^{n^2} c_1^n c_3 \frac{n(n+1)}{(n-1)!}  \frac{d}{\covol \Lambda}$$
where $c_1$ and $c_3$ are as in Theorem~{\rm{\TTerrorterms}}.
\endproclaim

To prove this proposition we will use the technique found, for example,
in the proof of \cite{\refMarcus, Lemma 2, p.~165}.
The technique requires that we express the boundary of $V_n$ as the
union of the images of easily-understood regions of $\BR^{n-1}$ under 
Lipschitz maps for which we have explicit Lipschitz factors;
by ``easily-understood'', we mean that we should be able
to find good approximations for the sizes of the intersections of the
regions with cubical lattices.
Thus, our proof will boil down to the two lemmas presented below, whose
proofs we will postpone until after the proof of Proposition~\PPlatticepoints.

Let $I_m'$ denote the subset
$\bigl\{\,(x_1,\ldots,x_m) \bigm| -2 \le x_1 \le x_2 \le \cdots \le x_m < 2\,\bigr\}$
of $\BR^m$, so that $I_m'$ can be obtained from the
simplex $I_m$ defined in Section~\SSregion\ by removing one face.

\proclaim{Lemma~\LLlatticepointsI}
Let $n>1$ be an integer, let $\delta$ be a real number such that
$\delta \le (n-1) 6^{-n}/c_6$ where $c_6 = \sqrt{3}/8$, and let
$M$ denote the shifted cubical lattice in $\BR^{n-1}$ with
center at $(-2,\ldots,-2)$ and edge length $\delta$.  Then
$$\#(M\cap I_{n-1}') \le c_7 \frac{4^{n-1}}{(n-1)!\,\delta^{n-1}},$$
where $c_7 = (1 + \sqrt{3}/162)^3$.
\endproclaim

The proof of this lemma appears immediately after the proof of 
Proposition~\PPlatticepoints.

Now we will define $n+1$ maps from $I_{n-1}$ to $\partial I_n$
whose images cover $\partial I_n$.
For every $i=0$, \dots,~$n$ let $\Delta_i\colon I_{n-1}\ra I_n$
be defined by 
$$\Delta_i(x_1,\ldots,x_{n-1}) = \cases
(-2,x_1,\ldots,x_{n-1})                                   &\text{if $i = 0$};\\
(x_1, \ldots, x_{i-1}, x_i, x_i, x_{i+1},\ldots, x_{n-1}) &\text{if $0 < i < n$;}\\
(x_1, \ldots, x_{n-1}, 2)                                 &\text{if $i = n$.}
\endcases$$
For every $m$ let $\| \cdot\| $ denote the sup-norm on $\BR^m$.

\proclaim{Lemma~\LLlipshitz}
For every index $i=0$, \dots,~$n$ and for every $\bx$ and $\By$ in $I_{n-1}$ we have
$\bigl\| (\Phi\circ\Delta_i)(\bx) - (\Phi\circ\Delta_i)(\By)\bigr\|
   \le c_6 6^n \| \bx - \By\| $,
where $c_6$ is as in Lemma~{\rm \LLlatticepointsI}.
\endproclaim

We will prove this lemma in the next section.

\demo{Proof of Proposition {\rm \PPlatticepoints}}
It is easy to verify
% checked 1/31/98
the statement of the proposition when $n=1$,
so henceforth we will assume that $n>1$.

Suppose $\Lambda$ is generated by the vectors $d_1 \be_1$, \dots,~$d_n \be_n$,
where the $\be_i$ are the standard unit vectors in $\BR^n$ and where the
$d_i$ are positive.
To every lattice point $\ell\in \Lambda$ we associate the ``brick''
$B_\ell = \ell + \bigl\{\,(x_1,\ldots,x_n) \bigm| \forall i: 0\le x_i<d_i\,\bigr\}$.
Let $S$ denote the set of all $\ell$ such that $B_\ell\subseteq V_n$, let
$T$ denote the set of all $\ell$ such that $B_\ell\cap V_n\neq\emptyset$, and
let $U = T\setminus S$.  It is clear that $\#T\ge \volume V_n / \covol \Lambda$
and that $\#S\le \volume V_n / \covol \Lambda$.  Since $\ell\in V_n$ implies
$\ell \in T$ we find that
$$\#(\Lambda\cap V_n) \le \#T = \#S + \#U \le \frac{\volume V_n}{\covol \Lambda} + \#U,$$
and since $\ell \in S$ implies $\ell \in V_n$ we find that
$$\#(\Lambda\cap V_n) \ge \#S = \#T - \#U \ge \frac{\volume V_n}{\covol \Lambda} - \#U,$$
so that
$$\left| \#(\Lambda\cap V_n) - \frac{\volume V_n}{\covol \Lambda}\right| \le \#U.$$
Thus, to prove the proposition it will be enough to show that 
$$\#U \le 6^{n^2} c_1^n c_3 \frac{n(n+1)}{(n-1)!}  \frac{d}{\covol \Lambda},$$
where $d$ is the mesh of $\Lambda$.

Let $\delta = (n-1) 6^{-n} d /c_6$, where $c_6$ is as in Lemma~\LLlatticepointsI,
and let $M$ be the shifted cubical lattice in $\BR^{n-1}$ with center
at $(-2,\ldots,-2)$ and with edge length $\delta$.  Notice that if $\bx$ is
a vector in $I_{n-1}$ then there is an element $\bm$ of $M\cap I_{n-1}'$
such that $\|\bx - \bm\| \le \delta$.

Suppose $\ell\in U$.  Then there is a point $\bz\in\partial V_n$
such that $\|\ell - \bz\|  \le d$.  Let $\By = \Phi^{-1}(\bz)$,
so that $\By\in\partial I_n$.  Finally, choose an $i$ and an $\bx\in I_{n-1}$
such that $\By = \Delta_i(\bx)$, and let $\bm$ be an element of $M\cap I_{n-1}'$
such that $\|\bx - \bm\| \le \delta$.
Using Lemma~\LLlipshitz, we find that
$$\align
\bigl\|\ell - (\Phi\circ\Delta_i)(\bm)\bigr\|
  & \le \|\ell - \bz\|  + \bigl\| \bz - (\Phi\circ\Delta_i)(\bm)\bigr\| \\
  & \le d + \bigl\| (\Phi\circ\Delta_i)(\bx) -  (\Phi\circ\Delta_i)(\bm)\bigr\| \\
  & \le d + c_6 6^n \delta\\
  & = nd.
\endalign
$$

We find that the number of elements of $U$ is bounded by the number of 
elements $\ell$ of $\Lambda$ for which there is an $i$ and an $\bm\in M\cap I_{n-1}'$
such that $\bigl\| \ell - (\Phi\circ\Delta_i)(\bm)\bigr\|  \le nd$.
Now, if $\bw\in\BR^n$ then the number of $\ell\in \Lambda$ 
with $\|\ell - \bw\|  \le nd$ is at most
$$\align
\left(\frac{2nd}{d_1}+1\right)\cdots\left(\frac{2nd}{d_n}+1\right)
 & = \frac{(2nd)^n}{\covol\Lambda} 
    \left(1+\frac{d_1}{2nd}\right)\cdots\left(1+\frac{d_n}{2nd}\right)\\
 & \le \frac{(2nd)^n}{\covol\Lambda} 
    \left(1+\frac{1}{2n}\right)^n \\
 & < \frac{(2nd)^n \exp(1/2)}{\covol\Lambda}.
\endalign
$$
Thus, 
using Lemma~\LLlatticepointsI\ we find that
$$\align
\#U 
 & \le (\#\text{\rm\ of possible $i$})\cdot(\#\text{\rm\ of possible $\bm$})
    \cdot(\#\text{\rm\ of $\ell$ for a given $i$ and $\bm$})\\
 &\le (n+1) \cdot\#(M\cap I_{n-1}')
     \cdot \frac{(2nd)^n \exp(1/2)}{\covol\Lambda}\cr
 &\le (n+1) c_7 \frac{4^{n-1}}{(n-1)! \,\delta^{n-1}}
     \frac{(2nd)^n \exp(1/2)}{\covol\Lambda} \\
 &= (n+1) c_7 \frac{4^{n-1}}{(n-1)!} 
    \left(\frac{c_6 6^n}{(n-1)d}\right)^{n-1}
     \frac{(2nd)^n \exp(1/2)}{\covol\Lambda}.
\endalign$$
Regrouping the terms in this last expression and
using the fact that $(n/(n-1))^{n-1}$ is less than $\exp(1)$, we find
$$\align
\#U &\le \frac{n (n+1)}{(n-1)!}
    6^{n^2}  \left(\frac{4 c_6}{3}\right)^n 
    \left(\frac{n}{n-1}\right)^{n-1}
    \left(\frac{c_7\exp(1/2)}{4 c_6}\right)
    \frac{d}{\covol\Lambda} \\
 &< \frac{n (n+1)}{(n-1)!}
    6^{n^2}  \left(\frac{4 c_6}{3}\right)^n 
    \left(\frac{c_7 \exp(3/2) }{4 c_6}\right)
    \frac{d}{\covol\Lambda} \\
 &= \frac{n(n+1)}{(n-1)!}
    6^{n^2}  c_1^n c_3
    \frac{d}{\covol\Lambda}.
\endalign$$
\myqed
\enddemo

We are left with the task of proving Lemma~\LLlatticepointsI.
Our proof will depend on the following simple fact.

\proclaim{Lemma \LLexactlattice}
Let $m>0$ be an integer, let $C\subset \BR^n$ be the cubical lattice
with edge-length $1/m$, and let $J_n'$ be the region
$\bigl\{\, (x_1,\ldots,x_n) \bigm| 0\le x_1\le \cdots \le x_n < 1\,\bigr\}$.
Then $\#(C\cap J_n') = \binom{n+m-1}{n}$.
\endproclaim

\demo{Proof}
The map $(y_1,\ldots,y_n) \mapsto (my_1+1, my_2+2, \ldots, my_n + n)$
gives a bijection between $C\cap J_n'$ and the set
$\bigl\{\, (z_1,\ldots,z_n) \in \BZ^n \bigm| 0 < z_1 < \cdots <  z_n < m+n\,\bigr\}$,
whose cardinality is equal to the number of ways one can choose
$n$ distinct integers from the set $\{\,1,2,\ldots,m+n-1\,\}$.
\myqed
\enddemo

\demo{Proof of Lemma~{\rm \LLlatticepointsI}}
The affine map $\bx \mapsto \bx/4 + (1/2,\ldots,1/2)$ sends $M$ to the
cubical lattice $C'\subset \BR^{n-1}$ with edge-length $\delta/4$ and
sends $I_{n-1}'$ to $J_{n-1}'$, so \hbox{$\#(M\cap I_{n-1}') = \#(C'\cap J_{n-1}')$.}
Let $m=\lceil 4/\delta \rceil$ and let $C\subset\BR^{n-1}$ be the 
cubical lattice with edge length $1/m$.
Then $1/m \le \delta/4$ and $m \le (4/\delta) (1+\delta/4)$, so we have
$$\align
\#(C'\cap J_{n-1}')
 & \le \#(C\cap J_{n-1}')\\
 & = \binom{m+n-2}{n-1}\\
 & = \frac{m^{n-1}}{(n-1)!} 
     \left(1 + \frac{1}{m}\right)\cdots\left(1 + \frac{n-2}{m}\right) \\
 & \le \frac{4^{n-1}}{(n-1)!\,\delta^{n-1}} \left(1+\frac{\delta}{4}\right)^{n-1}
   \left(1 + \frac{\delta}{4}\right)\cdots\left(1 + \frac{(n-2)\delta}{4}\right).
\endalign
$$
To finish our proof we need only show that 
$$
\left(1+\frac{\delta}{4}\right)^{n-1}
   \left(1 + \frac{\delta}{4}\right)\cdots\left(1 + \frac{(n-2)\delta}{4}\right)
\le c_7
\eqno{\EQmessyinequality}$$
when $n>1$ and $\delta\le (n-1)6^{-n}/c_6$.
Inequality~\EQmessyinequality\  is easy to verify when $n=2$
%checked on 1/31/98
and $n=3$,
%checked on 1/31/98
so we are left with the case $n\ge4$.

We have 
$$\align
\left(1+\frac{\delta}{4}\right)^{n-1}
   \left(1 + \frac{\delta}{4}\right)\cdots\left(1 + \frac{(n-2)\delta}{4}\right)
& \le \left(1+\frac{\delta}{4}\right)^{n-1}
  \left(1 + \frac{(n-2)\delta}{4}\right)^{n-2} \\
& <\exp\bigl((n-1)\delta/4\bigr) \exp\bigl((n-2)^2\delta/4\bigr) \\
& = \exp\left( (n^2-3n+3)\delta/4\right)\\
& \le \exp\left( \frac{2(n-1)(n^2-3n+3)}{6^n\sqrt{3}}\right),
\endalign
$$
and it is easy to verify that this last expression is less than $c_7$ when $n\ge 4$.
%checked on 1/31/98
Thus inequality~\EQmessyinequality\ holds for all $n>1$ and we are done.
\myqed
\enddemo

%%%%%%%%%%%%%%%%%%%%%%%%%%%%%%%%%%%%%%%%%%%%%%%%%%%%%%%%%%%%%%%%%%%%%%%%
\subhead\SSlipshitz. Proof of the Lipschitz bound
\endsubhead

Since $\Phi = \Chi\circ\Psi$, 
Lemma~\LLlipshitz\  follows from the following more precise result.

\proclaim{Lemma \LLseverallipshitz}
For every $n>1$ the following statements hold{\rm{:}}
\roster
\item"1."
For every index $i=0$, \dots,~$n$ and for every $\bx$ and $\By$ in $I_{n-1}$ we have
${\bigl\|  \Delta_i(\bx)- \Delta_i(\By) \bigr\| } = {\|\bx - \By\| }$.
\item"2."
For every $\bx$ and $\By$ in $I_n$ we have
$\bigl\|  \Psi(\bx)- \Psi(\By) \bigr\|  
\le 3^{n-1} \sqrt{n} \ \|  \bx - \By \| $.
\item"3."
For every $\bx$ and $\By$ in $\BR^n$ we have
$\bigl\|  \Chi(\bx)- \Chi(\By) \bigr\| 
\le \bigl(3\sqrt{3}\cdot2^{n-3}/\sqrt{n}\bigr) \|  \bx - \By \| $.
\endroster
\endproclaim

\demo{Proof}
The first statement of the lemma follows immediately from the definitions of $\Delta_i$ and
of the sup-norm.  To prove the second and third statements, we will use
two basic facts.

{\it First fact}:  Suppose $R$ is a convex open region of $\BR^m$ and $f$ is
a continuous function from the closure of $R$ to $\BR^n$ that is differentiable 
on $R$.  
Suppose $M$ is a real number
such that for all $\bz\in R$ and for all $\bx\in \BR^m$ we have
$\bigl\| D\cdot\bx \bigr\|\le M \|\bx\|$,
where $D$ is the Jacobian matrix of $f$ evaluated at $\bz$.
Then for all $\bx$ and $\By$ in the closure of $R$ we have 
$\bigl\| f(\bx) - f(\By)\bigr\| \le M \|\bx-\By\|$.
(This statement is an immediate consequence of \cite{\refApostol, Theorem~12.9, p.~355}.)

{\it Second fact}: Suppose $D = [d_{i,j}]$ is an $n\times m$ matrix of real numbers,
and let 
$$M = \max_{1\le i\le n} \ \sum_{j=1}^m |d_{i,j}|.$$
Then for all $\bx\in \BR^m$ we have
$\bigl\| D\cdot\bx \bigr\| \le M \| \bx\|$.
(We leave the simple verification of this fact to the reader.)

Let us proceed to the proof of the second statement of the lemma.
Suppose $\br = (r_1,\ldots, r_n)$ is an element of $I_n$, and let $D = [d_{i,j}]$
be the Jacobian of $\Psi$ evaluated at $\br$.
The definition of $\Psi$ shows that $d_{i,j}$ is $(-1)^i$ times
the $(i-1)$th symmetric polynomial in the
$n-1$ numbers $\{\,r_k \mid k\neq j\,\}$, and since $\br\in I_n$ the
absolute value of this entry is at most $\binom{n-1}{i-1} 2^{i-1}$.
Thus 
$$\max_{1\le i\le n} \ \sum_{j=1}^n |d_{i,j}|
\le \max_{1\le i\le n} \ n\binom{n-1}{i-1}2^{i-1}.$$
Let $b_i= n\binom{n-1}{i-1}2^{i-1}$.  Then for $i>1$ we have
$b_i/b_{i-1} = 2 (n-i+1)/(i-1)$, so $b_i > b_{i-1}$ if and only if
$i\le 2n/3 + 1$.  We see that the maximum value of $b_i$ occurs
when $i = \lfloor 2n/3 \rfloor + 1$.  Thus Lemma~\LLtwothirds\ below
shows that $b_i \le 3^{n-1}\sqrt{n}$ for all~$i$, and this result,
combined with first and second facts above, proves
the second statement of the lemma.

We are left to prove the third statement.
It is not hard to see from the definition of $\Chi$ that $\Chi$ is
an affine map; that is, there is an $n\times n$ matrix $D$ such that
$\Chi(\bc) = \Chi\bigl( (0,\ldots,0)\bigr) + D\cdot\bc$ for every $\bc \in \BR^n$.
Furthermore, it is easy to check that  every entry of $D$ is non-negative.
It follows that the maximum row-sum of $D$ is equal to the
sup-norm of $\Chi\bigl( (1,\ldots,1)\bigr) - \Chi\bigl((0,\ldots,0)\bigr)$.
Let $g_1 = x^n\big( (x+1/x)^n + \cdots + (x+1/x) + 1\big)$
and let $g_0 = x^n (x+1/x)^n$.
The definition of $\Chi$ 
in terms of polynomials
shows that
the sup-norm of $\Chi\bigl( (1,\ldots,1)\bigr) - \Chi\bigl((0,\ldots,0)\bigr)$
is equal to the largest coefficient of the polynomial $g_1 - g_0$,
and it is not hard to check that this largest coefficient is
the coefficient of $x^n$.  In other words,
the maximum row-sum of $D$ is equal to the sum over the even integers $j$ less than
$n$ of $\binom{j}{j/2}$.
Lemma~\LLonehalf\ below shows that this sum is at most 
$3\sqrt{3}\cdot 2^{n-3}/\sqrt{n}$.
Using the first and second facts above, we see that the third statement of the lemma
is true.
\myqed
\enddemo

\proclaim{Lemma \LLtwothirds}
Let $n$ be a positive integer and let $i = \lfloor 2n/3 \rfloor + 1$.
Then we have
$n\binom{n-1}{i-1} 2^{i-1} \le 3^{n-1}\sqrt{n}$.
\endproclaim

\demo{Proof}
For every $n>0$ let
$r_n = \bigl(n \binom{n-1}{i-1} 2^{i-1}\bigr)\big/\bigl(3^{n-1}\sqrt{n}\bigr)$,
where $i = \lfloor 2n/3 \rfloor+1$.
Then we find that for every $m\ge 0 $ we have
$$\frac{r_{3m+4}}{r_{3m+1}} = \frac{(m+2/3)\sqrt{(m+1/3)(m+4/3)}}{(m+1/2)(m+1)} < 1$$
and
$$\frac{r_{3m+3}}{r_{3m+1}} = \frac{(m+2/3)\sqrt{m+1/3}}{(m+1/2)\sqrt{m+1}} < 1$$
and 
$$\frac{r_{3m+2}}{r_{3m+1}} = \frac{\sqrt{(m+1/3)(m+2/3)}}{m+1/2} < 1.$$
%checked 1/31/98
The first inequality shows that $r_n\le r_1$ for all $n\equiv 1\pmod3$,
and then the second and third inequalities show that $r_n\le r_1$ for all $n$.
Since $r_1 = 1$, the lemma follows.
\myqed
\enddemo

\proclaim{Lemma \LLonehalf}
We have the following inequalities{\rm{:}}
\roster
\item"1."
If $j$ is a positive integer, then
$\binom{2j}j < 4^j/\sqrt{\pi j}$.
\item"2."
Let $n$ be a positive integer and let $i = \lfloor (n-1)/2 \rfloor$.
Then 
$$\binom00 + \binom21 + \cdots + \binom{2i}i 
 \le \frac{3\sqrt{3}}{8} \frac{2^n}{\sqrt{n}}.$$
\endroster
\endproclaim

\demo{Proof}
For every $j>0$ let $r_j = \binom{2j}j \sqrt{j}/4^j$.  Then 
$$\frac{r_{j+1}}{r_j} = \frac{j+1/2}{\sqrt{j(j+1)}} > 1$$
so the $r_j$ form an increasing sequence.  Stirling's formula
shows that the $r_j$ approach $1/\sqrt{\pi}$, and the first inequality
of the lemma follows.

Let $c_8 = 9\sqrt{2}/16$.  To prove the second part of the lemma, we will first show that
for every $m>0$ we have
$$\binom00 + \binom21 + \cdots + \binom{2m}m \le c_8 \frac{4^m}{\sqrt{m}}. \eqno{\EQbinsum}$$
One checks by hand that~\EQbinsum\ holds for $m=1$, $m=2$, and $m=3$.
%checked 1/31/98
Let $s_m$ denote the sum on the left hand side of inequality~\EQbinsum,
and suppose that~\EQbinsum\ holds when $m$ is equal to some integer $j>2$.
Then by using the first part of the lemma we find that
$$\align
s_{j+1} 
 &  =  s_j + \binom{2j+2}{j+1} \\
 & \le c_8 \frac{4^j}{\sqrt{j}} + \frac1{\sqrt{\pi}}\frac{4^{j+1}}{\sqrt{j+1}} \\
 &  =  c_8 \frac{4^{j+1}}{\sqrt{j+1}}\left(\frac{\sqrt{j+1}}{4\sqrt{j}}+ \frac1{c_8\sqrt{\pi}}\right)\\
 &  <  c_8 \frac{4^{j+1}}{\sqrt{j+1}},
\endalign
$$
where the final inequality holds because the expression in parentheses is less than~$1$
when $j\ge 3$.  
%checked on 2/2/98 -- the expression is just under 0.998.
By induction, inequality~\EQbinsum\ holds for all $m>0$.

Now suppose that $n>4$ is an integer, and write $n = 2m+r$ where $m>1$ is an integer
and $r=1$ or $r=2$.
Then the $i$ in the second statement of the lemma is equal to $m$,
and using inequality~\EQbinsum\ we see that
$$\binom00 + \binom21 + \cdots + \binom{2i}i 
 \le c_8 \frac{4^{m}}{\sqrt{m}}
 =  c_8 2^{-r}\sqrt{\frac{n}{m}}\frac{2^n}{\sqrt{n}}
 < \frac{3\sqrt{3}}{8} \frac{2^n}{\sqrt{n}},$$
where the last inequality follows from the fact that
$c_8 2^{-r}\sqrt{(2m+r)/m} < 3\sqrt{3}/8$
when $m>1$ and $r=1$ or $r=2$.
%checked 1/31/98
Thus the second inequality of the lemma holds for all $n>4$.
But direct computation shows that the inequality holds for $n \le 4$ as well.
%checked 1/31/98
\myqed
\enddemo

%%%%%%%%%%%%%%%%%%%%%%%%%%%%%%%%%%%%%%%%%%%%%%%%%%%%%%%%%%%%%%%%%%%%%%%%
\subhead\SSlowerbounds.  Nonzero lower bounds for lattice points
\endsubhead

If one is interested in obtaining nonzero lower bounds
on the size of the intersection of a rectilinear lattice with
one of the regions $V_n$  --- as we will be in Part~\PartAV\ ---
then Proposition~\PPlatticepoints\ is only helpful
when the mesh of the lattice is very small, say on the 
order of $6^{-n^2}$.
In this section we will prove several lemmas that can be
used to give nontrivial lower bounds even when the mesh
of the lattice is close to $1$.
The idea is to show that the complicated region $V_n$ contains
a simple ``diamond-shaped'' region, and to give a good
lower bound on the size of the intersection of a rectilinear lattice with such 
a diamond-shaped region.
We will demonstrate the use of these lemmas in the proof of Theorem~\TTlowerbounds\ in
Section~\SSisogenyclass.  Lemma~\LLniceregion\ will also be the critical element in the
proof of Theorem~\TTgrouporders.

\proclaim{Lemma~\LLniceregion}
Suppose $\bb = (b_1,\ldots, b_n) \in \BR^n$ satisfies
$| b_n/2 | + \sum_{i=1}^{n-1} |b_i| \le 1$.
Then $\bb\in V_n$.
\endproclaim

\demo{Proof}
Let $W_n$ be the set of all vectors 
that satisfy the hypothesis of the lemma, and let $W_n^{\roman{o}}$ denote the interior of $W_n$.
Since $V_n$ is closed, to prove that $W_n\subseteq V_n$ it will be enough to
prove that $W_n^{\roman{o}}\subseteq V_n$, and to accomplish this it will be enough
to show that $W_n^{\roman{o}}\cap V_n\neq\emptyset$
and that $W_n^{\roman{o}}\cap\partial V_n = \emptyset$.  The first statement is clear,
because the vector $(0,\ldots,0)$ lies in both $W_n^{\roman{o}}$ and $V_n$.
To prove the second statement, we use Lemma~\LLboundary\ as follows.

Suppose that $\bb=(b_1,\ldots,b_n)$ is an element of $W_n^{\roman{o}}\cap V_n$,
and suppose $z$ is a root of $g_\bb$.  We calculate that
$$\align
g'_\bb(z) & = 2nz^{2n-1} + b_1\bigl((2n-1)z^{2n-2} + 1\bigr) + b_2\bigl((2n-2)z^{2n-3}+2z\bigr) + \cdots\\
          & \qquad\qquad  + b_{n-1}\bigl((n+1)z^n + (n-1)z^{n-2}\bigr) + b_n n z^{n-1},
\endalign
$$
and since all the roots of $g_\bb$ lie on the unit circle it follows that
$$\align
\bigl|g'_\bb(z)\bigr| & \ge 2n - |b_1|\bigl((2n-1)+1\bigr) - |b_2|\bigl((2n-2)+2\bigr) - \cdots \\
            & \qquad\qquad - |b_{n-1}|\bigl((n+1)+(n-1)\bigr) - |b_n| n\\
            & = 2n \left( 1 -\left|\frac{b_n}{2}\right| - \sum_{i=1}^{n-1} |b_i| \right)\\
            & > 0
\endalign
$$
because $\bb \in W_n^{\roman{o}}$.  Thus $g_\bb$ has no multiple roots, so 
$\bb$ cannot be in $\partial V_n$, and the lemma is proved.
\myqed
\enddemo

\proclaim{Lemma \LLwedge}
Let $\Lambda\subset\BR^n$ be a rectilinear lattice, let $r>0$ be a real
number, and let $W\subset\BR^n$ be the region
$W=\bigl\{\,(x_1,\ldots,x_n) \bigm| x_1+\cdots + x_n \le r \quad \&\quad  \forall i: x_i\ge 0\,\bigr\}$.
Then $$\#(\Lambda\cap W)\ge \frac{r^n}{n!\, \covol\Lambda}.$$
\endproclaim

\demo{Proof}
Recall that to every lattice point $\ell\in\Lambda$ we associate a brick $B_\ell$.
It is easy to see that for every $x\in W$ there is an $\ell\in\Lambda\cap W$ such that
$x\in B_\ell$. Thus $\#(\Lambda\cap W)$ is at least the ratio of the volume of $W$
to the volume of a brick.  The lemma then follows from the fact that the volume of $W$
is $r^n/n!$.
\myqed
\enddemo

Let $r>0$ be a real number, and let
$U=\bigl\{\,(x_1,\ldots,x_n) \bigm| |x_1|+\cdots + |x_n| \le r\,\bigr\}$.
Suppose $\Lambda$ is a rectilinear lattice
generated by the vectors $d_1 \be_1, \ldots, d_n \be_n$,
where the $\be_i$ are the standard unit vectors in $\BR^n$ and where the
$d_i$ are positive.
For every subset $S$ of $\{\,1,\ldots,n\,\}$ let $d_S$ denote the sum $\sum_{i\in S} d_i$.

\proclaim{Lemma \LLdiamond}
We have
$$\#(\Lambda\cap U) \ge \frac{1}{n!\,\covol\Lambda}
\sum_{S\colon d_S \le r} (r - d_S)^n.$$
\endproclaim

\demo{Proof}
For every subset $S$ of $\{\,1,2,\ldots,n\,\}$ let 
$U_S$ denote the subset of $U$ consisting of those
$(x_1,\ldots,x_n)$ such that $x_i \ge 0$ if $i\not\in S$ and $x_i\le -d_i$ if $i\in S$.
Furthermore, let 
$$W_S=\bigl\{\,(y_1,\ldots,y_n) \bigm| y_1+\cdots + y_n \le r-d_S \quad \&\quad  \forall i: y_i\ge 0\,\bigr\}.$$
Clearly $\Lambda\cap U$ is the disjoint union of the $\Lambda\cap U_S$.  Also,
for every $S$ there is a bijection between $\Lambda\cap U_S$ and $\Lambda\cap W_S$
given by the map $(x_1,\ldots,x_n)\mapsto(y_1,\ldots,y_n)$, 
where $y_i = x_i$ if $x_i\ge 0 $ and $y_i = -x_i - d_i$ if $x_i<0$.
The lemma then follows from Lemma~\LLwedge.
\myqed
\enddemo

%%%%%%%%%%%%%%%%%%%%%%%%%%%%%%%%%%%%%%%%%%%%%%%%%%%%%%%%%%%%%%%%%%%%%%%%
\head\PartAV. Abelian varieties over finite fields
\endhead

%%%%%%%%%%%%%%%%%%%%%%%%%%%%%%%%%%%%%%%%%%%%%%%%%%%%%%%%%%%%%%%%%%%%%%%%
\subhead\SSweilpolys.  Properties of Weil polynomials
\endsubhead

To every abelian variety $A$ over $\BF_q$
one associates its characteristic polynomial of Frobenius $f_A\in\BZ[x]$, sometimes called
the {\it Weil polynomial\/} or the {\it Weil $q$-polynomial\/}
of the variety.  The polynomial $f_A$ is monic of degree
twice the dimension of $A$, and Weil's ``Riemann Hypothesis''
says that all of its roots in $\BC$ have magnitude $\sqrt{q}$.
Furthermore, the Honda-Tate theorem (see \cite{\refTate})
implies that the real roots of $f_A$, if there are any, have even multiplicity.
It follows that $f_A$ can be written
$$f_A=\bigl(x^{2n} + q^n\bigr) + a_1 \bigl(x^{2n-1}+q^{n-1}x\bigr) + \cdots 
         + a_{n-1}\bigl(x^{n+1}+qx^{n-1}\bigr) + a_n x^n$$
for some integers $a_i$. 
The variety $A$ is called {\it ordinary}, and $f_A$ is called an 
{\it ordinary Weil $q$-polynomial}, if the middle coefficient $a_n$ is
coprime to $q$.

In this section we will prove two propositions, one to let us
easily identify ordinary Weil polynomials and the other to provide a
necessary condition for a polynomial to be a non-ordinary Weil polynomial.
We will use properties of Newton polygons in our proofs; see 
for example \cite{\refGoss, Chapter~2} for the basic facts about Newton
polygons that we will require.

\proclaim{Proposition \PPordinary}
Suppose $f\in\BZ[x]$ is a monic polynomial of degree $2n$ all of whose roots
in $\BC$ have magnitude $\sqrt{q}$, and suppose the middle coefficient of
$f$ is coprime to $q$. Then $f$ is an ordinary Weil $q$-polynomial.
\endproclaim

\demo{Proof}
It is easy to check that every irreducible factor of $f$ in $\BZ[x]$ must have
even degree and have middle coefficient coprime to $q$, and if these factors
are ordinary Weil polynomials then so must be $f$.  Thus it suffices to consider
the case where $f$ is irreducible. 
The only real numbers with magnitude $\sqrt{q}$ are $\sqrt{q}$ and $-\sqrt{q}$,
and the minimal polynomials over $\BQ$ of these numbers do not satisfy the hypotheses
of the proposition, so we reduce the proof to the
case where $f$ is irreducible and has no real roots.

Let $v$ be the valuation on $\BQpb$ such that $v(q)=1$.
The Honda-Tate theorem says that there is a unique integer $e$ such that $f^e$ is
the Weil polynomial of a simple abelian variety, and the theorem
tells us how to calculate $e$:  It is the smallest positive integer 
such that for every root $\pi$ of $f$ in $\BQpb$ the rational number
$e v(\pi) [\BQp(\pi):\BQp]$ is an integer.  (If $f$ had real roots, $e$ would
be the smallest {\it even\/} integer satisfying this condition.)
But the slopes of the Newton polygon (with respect to $v$) for our $f$ are $0$ and $-1$, 
so $v(\pi)$ is either $0$ or $1$, and it follows that $e=1$.
Thus $f$ is an ordinary Weil polynomial.
\myqed
\enddemo

\proclaim{Lemma \LLnewton}
Let $q$ be a power of a prime $p$,
let $v$ be the $p$-adic valuation on $\BQpb$ normalized so that $v(q)=1$,
and let $f$ be a Weil $q$-polynomial.  Then the vertices of the Newton polygon
for $f$ {\rm{(}}with respect to $v${\rm{)}} are integer lattice points.
\endproclaim

\demo{Proof}
If the Newton polygons of two polynomials have vertices that
are integer lattice points, then so does the Newton polygon of their product.
Thus we need only prove the lemma when $f$ is the Weil polynomial of a simple
abelian variety.  In this case $f=P^e$ for some irreducible polynomial $P$.
Suppose $P$ factors in $\BZ_p[x]$ into the product of irreducibles $P_i$.
As we mentioned in the preceding proof, the Honda-Tate theorem shows that
if $\pi_i$ is a root of $P_i$, then $e v(\pi_i) \deg P_i$ is an integer.
But the Newton polygon of $P_i^e$ is a straight line from $(0, e v(\pi_i) \deg P_i)$
to $(\deg P_i, 0)$, so its vertices are integer lattice points.  Thus the
vertices of the Newton polygon for $f=P^e$ are integer lattice points as well.
\myqed
\enddemo

\proclaim{Proposition \PPnonordinary}
Suppose $A$ is a non-ordinary $n$-dimensional abelian variety over $\BF_q$,
let $a_n$ be the middle coefficient of its Weil polynomial, and let 
$v$ be the $p$-adic valuation on $\BQ$ with $v(q) = 1$.  Then $v(a_n) \ge 1/2$.
\endproclaim

\demo{Proof}
Since $f_A$ has the form
$$f_A=\bigl(x^{2n} + q^n\bigr) + a_1 \bigl(x^{2n-1}+q^{n-1}x\bigr) + \cdots 
         + a_{n-1}\bigl(x^{n+1}+qx^{n-1}\bigr) + a_n x^n$$
we see that if $(n+i,j)$ is a vertex of the Newton polygon for $f_A$, with $i>0$,
then so is $(n-i, j+i)$.  

Let $i$ be the smallest non-negative integer for which there is a $j$ such that
$(n+i,j)$ is a vertex of the Newton polygon for $f_A$.  If $i=0$ then $v(a_n) = j$
is an integer, and since $A$ is not ordinary $v(a_n)\ge 1>1/2$.  On the other
hand, if $i>0$ then
the point $(n,v(a_n))$ lies above the line connecting $(n-i,j+i)$ to $(n+i,j)$,
and it follows that $v(a_n)\ge j+i/2 \ge 1/2$.
\myqed
\enddemo

%%%%%%%%%%%%%%%%%%%%%%%%%%%%%%%%%%%%%%%%%%%%%%%%%%%%%%%%%%%%%%%%%%%%%%%%
\subhead\SSisogenyclass.  Counting isogeny classes
\endsubhead

In this section we will prove Theorems~\TTerrorterms\ and~\TTlowerbounds\ from 
the introduction.  We begin with Theorem~\TTerrorterms.
Let $q$ be a power of a prime $p$ and let $n$ be a positive integer.
Note that the conclusions of Theorem~\TTerrorterms\ follow easily when $n=1$,
because in that case $\#\CO(q,n)$ is simply the number of integers
$t$ such that $|t| \le 2\sqrt{q}$ and $(t,q)=1$, and
$1\le\#\CN(q,n)\le 5$ (as follows from \cite{\refWaterhouse, Theorem~4.1}, for example).
So henceforth we will assume that $n>1$.
Also, the statement that $\#\CO(q,n-1)\le\#\CN(q,n)$
follows from the existence of 
the injection from $\CO(q,n-1)$ to $\CN(q,n)$ obtained by
multiplying an ordinary isogeny class by the isogeny class of a fixed
supersingular elliptic curve.  We are left to prove the other two
inequalities of the theorem.

Let $\be_1$, \dots,~$\be_n$ denote the standard basis vectors of $\BR^n$.
Our arguments will involve three lattices in $\BR^n$: The first lattice,
denoted $\Lambda_q$, is generated by the vectors $q^{-i/2}\be_i$; the second,
denoted $\Lambda'_q$, is generated by the same set of vectors, except with
$q^{-n/2}\be_n$ replaced with $pq^{-n/2}\be_n$;  and the third lattice,
denoted $\Lambda''_q$, is generated by the same set as was $\Lambda_q$, but with
$q^{-n/2}\be_n$ replaced with $sq^{-n/2}\be_n$, where $s$ is the smallest
power of $p$ such that $q \mid s^2$.  
Thus $\Lambda_q\supset\Lambda'_q\supseteq\Lambda''_q$.

We noted earlier that if $A$ is an $n$-dimensional abelian variety
over $\BF_q$ then its Weil polynomial $f_A$ has all complex roots on the
circle $|z| = \sqrt{q}$ and its real roots have even multiplicity.
If we write
$$f_A=\bigl(x^{2n} + q^n\bigr) + a_1 \bigl(x^{2n-1}+q^{n-1}x\bigr) + \cdots 
       + a_{n-1}\bigl(x^{n+1}+qx^{n-1}\bigr) + a_n x^n$$
and let 
$\bb = \bigl(a_1 q^{-1/2}, a_2 q^{-1}, \ldots, a_n q^{-n/2}\bigr)$,
then $\bb\in\Lambda_q$ and in the notation of Section~\SSregion\ we
have $f_A(x) = q^n g_\bb(x/\sqrt{q})$.
Furthermore, $g_\bb$ has all of its roots on the unit circle, and its real roots
have even multiplicity, so $\bb\in V_n$.  

The Honda-Tate theorem shows that
this association $A\mapsto\bb$ gives us an injection $\Theta$ from the set $\CI(q,n)$
to $\Lambda_q\cap V_n$.  Proposition~\PPordinary\ shows that $\Theta$
gives a bijection between $\CO(q,n)$ and 
$(\Lambda_q\cap V_n)\setminus(\Lambda'_q\cap V_n)$, and Proposition~\PPnonordinary\ shows
that $\Theta(\CN(q,n))\subseteq \Lambda''_q\cap V_n$.
Thus, Theorem~\TTerrorterms\ is a consequence of the following:

\proclaim{Proposition \PPAVlattice}
Let $q$ be a power of a prime $p$ and let $n>1$ be an integer.
Then we have
$$\bigl| \#(\Lambda_q \cap V_n) - \#(\Lambda'_q \cap V_n) - v_n r(q) q^{n(n+1)/4}\bigr|
 \le 6^{n^2} c_1^n c_2 \frac{n(n+1)}{(n-1)!} q^{(n+2)(n-1)/4}$$
and
$$ \#(\Lambda''_q\cap V_n) \le \left(v_n+6^{n^2} c_1^n c_3 \frac{n(n+1)}{(n-1)!}\right) q^{(n+2)(n-1)/4}.$$
\endproclaim

\demo{Proof}
Note that the lattice $\Lambda_q$ has covolume $q^{-n(n+1)/4}$ and mesh $q^{-1/2}$.  Also,
the lattice $\Lambda'_q$ has covolume $pq^{-n(n+1)/4}$, and its mesh is $q^{-1/2}$ unless
$n=2$ and $q=p$, in which case its mesh is $1$.  Applying Proposition~\PPlatticepoints\ to
these two lattices, combining the resulting inequalities, and using the
fact that $r(q) = 1 - 1/p$, we find that
the left-hand side of the first inequality of the proposition is at most 
$$6^{n^2} c_1^n c_3 \frac{n(n+1)}{(n-1)!} q^{(n+2)(n-1)/4}\left(1 + \frac{d\sqrt{q}}{p}\right)$$
where $d = 1$ if $n=2$ and $q=p$, and $d=q^{-1/2}$ otherwise.  We see
that $$1 + \frac{d\sqrt{q}}{p} \le 1 + \frac{1}{\sqrt{p}} \le 1 + \frac{1}{\sqrt{2}},$$ and since
$c_2 = c_3 (1 + 1/\sqrt{2})$ we obtain the first inequality of the proposition.

The lattice $\Lambda''_q$ has covolume $sq^{-n(n+1)/4}$ and its mesh is at most $1$,
so Proposition~\PPlatticepoints\ tells us that
$$ \#(\Lambda''_q\cap V_n) \le \left(v_n+6^{n^2} c_1^n c_3 \frac{n(n+1)}{(n-1)!}\right) \frac{q^{n(n+1)/4}}{s}.$$
But $s\ge q^{1/2}$, so we obtain the second 
inequality of the proposition.
\myqed
\enddemo

Now we turn to Theorem~\TTlowerbounds.  Again, the theorem is easy to prove when $n=1$, so we
will assume that $n>1$.
Let $I\subset \BR$ be the interval $[-1/n, 1/n]$ and let $U\subset\BR^{n-1}$ be the region
$$U=\left\{\,(x_1,\ldots,x_{n-1}) \biggm| |x_1|+\cdots + |x_{n-1}| \le 1 - \frac{1}{2n} \,\right\}.$$
Then Lemma~\LLniceregion\ shows that $U\times I\subset V_n$.  Let $\Lambda$ be the
lattice in $\BR^{n-1}$ generated by $q^{-1/2}\be_1, \ldots, q^{-(n-1)/2}\be_{n-1}$ and
let $M$ be the set $$M = \bigl\{\,mq^{-n/2} \bigm| m \in \BZ \quad\&\quad (m,q) = 1\,\bigr\}.$$
Then we have 
$$ (\Lambda\cap U)\times(M\cap I) \subseteq (\Lambda_q\cap V_n)\setminus(\Lambda'_q\cap V_n)$$
so the product of $\#(\Lambda\cap U)$ with $\#(M\cap I)$ gives a lower bound on
$\#\CO(q,n)$.
A simple argument shows that $\#(M\cap I) \ge 2(1 - 1/p) q^{n/2}/n - 2$,  so the following lemma
completes the proof of Theorem~\TTlowerbounds.

\proclaim{Lemma \LLAVdiamond}
We have
$$\#(\Lambda\cap U) >  (c_4/2) (c_5 n)^{-2(\log 2)/(\log q)} \frac{2^n}{(n-1)!}  q^{n(n-1)/4}.$$
\endproclaim

\demo{Proof}
Let $j$ be the smallest element of the set consisting of $n$ and
those integers greater than $2(\log c_5 n)/(\log q)$,
and let  $T = \{\, j, j+1, \ldots, n-1\,\}$, so that $T=\emptyset$ if $j = n$.
Suppose $S\subseteq T$.  If $j=n$ then $d_S = 0 < 1/n$ (in the notation of Lemma~\LLdiamond).
If $j<n$ then 
$$d_S\le d_T < \frac{q^{-j/2}}{1-1/\sqrt{q}}\le c_5 q^{-j/2} \le \frac{1}{n},$$
because $j\ge 2(\log c_5n)/(\log q)$.
Applying Lemma~\LLdiamond\ to $\Lambda$ and $U$ (so that $r = 1-1/2n$), we find that 
$$\align
\#(\Lambda\cap U) & \ge \frac{q^{n(n-1)/4}}{(n-1)!} \sum_{S\colon d_S \le r} (r - d_S)^{n-1} \\
& \ge \frac{q^{n(n-1)/4}}{(n-1)!} \sum_{S\subseteq T} \left( 1 - \frac{3}{2n}\right)^{n-1} \\
& \ge \frac{q^{n(n-1)/4}}{(n-1)!} 2^{n-j} \left( 1 - \frac{3}{2n}\right)^{n-1} \\
& > \frac{q^{n(n-1)/4}}{(n-1)!} 2^{n-j} c_4
\endalign
$$
because $(1-3/2n)^{n-1} > \exp(-3/2) = c_4$.
To complete the proof of the lemma we need only show that
$$2^{-j} \ge (1/2) (c_5 n)^{-2(\log 2)/(\log q)},$$
but this follows directly from the inequality $j < 1 + 2(\log c_5 n)/(\log q)$.
\myqed
\enddemo

%%%%%%%%%%%%%%%%%%%%%%%%%%%%%%%%%%%%%%%%%%%%%%%%%%%%%%%%%%%%%%%%%%%%%%%%
\subhead\SSgrouporders.  Group orders of abelian varieties
\endsubhead

In this section we will prove Theorem~\TTgrouporders.  We begin with a simple lemma.

\proclaim{Lemma~\LLWeilpolys}
Let $q$ be a power of a prime number and let $n>0$ be an integer.
Suppose $a_1$, \dots,~$a_n$ are integers such that 
$$\left|\frac{a_n}{2q^{n/2}}\right| + \sum_{i=1}^{n-1} \left|\frac{a_i}{q^{i/2}}\right| \le 1$$
and $(a_n,q)=1$.
Then
$$f=\bigl(x^{2n} + q^n\bigr) + a_1 \bigl(x^{2n-1}+q^{n-1}x\bigr) + \cdots 
     + a_{n-1}\bigl(x^{n+1}+qx^{n-1}\bigr) + a_n x^n$$
is an ordinary Weil $q$-polynomial.
\endproclaim

\demo{Proof}
For every $i$ let $b_i= a_i/q^{i/2}$ and let $\bb = (b_1,\ldots,b_n)$.
Then we have $f(x) = q^n g_\bb(x/\sqrt{q})$.  By Lemma~\LLniceregion\ every root of
$g_\bb$ lies on the unit circle, so every root of $f$ has magnitude $\sqrt{q}$.
Since $(a_n,q)=1$, Proposition~\PPordinary\ tells us that $f$ is an ordinary Weil polynomial.
\myqed
\enddemo

\demo{Proof of Theorem {\rm \TTgrouporders}}
If $A$ is an abelian variety over $\BF_q$ then $\#A(\BF_q) = f_A(1)$,
so to prove the theorem it will be enough to find an ordinary Weil $q$-polynomial
$f$ of degree $2n$ with $f(1)=m$.  We have slightly different arguments for
the cases $n>2$ and $n=2$;  let us begin by proving the theorem when $n>2$.

We will choose the coefficients $a_1$, \dots,~$a_n$ of $f$ one at a time. To begin, we 
let $g_1 = m - \bigl(q^n+1\bigr)$ and pick $a_1$ so that the absolute value of
$g_1 - a_1(q^{n-1}+1)$ is minimized.
Now suppose we have chosen $a_1$ through $a_{i-1}$, for some $i<n$.
Let $$g_i = m - \bigl(q^n+1\bigr) - a_1\bigl(q^{n-1}+1\bigr) - \cdots - a_{i-1}\bigl(q^{n-i+1}+1\bigr),$$
and pick $a_i$ so that the absolute value of $g_i - a_i(q^{n-i}+1)$ is minimized.
Finally, pick $a_n$ so that
$$m = \bigl(q^n+1\bigr) + a_1\bigl(q^{n-1}+1\bigr) + \cdots 
      + a_{n-1}\bigl(q+1\bigr) + a_n. \eqno{\EQnumpoints}$$

In a moment we may change the values of $a_{n-1}$ and $a_n$,
but let us first deduce some properties of the $a_i$ as they stand.
The hypotheses of the theorem show that 
$$\Bigl| m - \bigl(q^n + 1\bigr)\Bigr| < C_q\sqrt{q}\bigl(q^{n-1}+1\bigr)
= \Bigl(\bigl\lfloor B_q\sqrt{q}\bigr\rfloor + 1/2\Bigr) \bigl(q^{n-1}+1\bigr),$$
so we have $$|a_1|\le \bigl\lfloor B_q\sqrt{q}\bigr\rfloor \le B_q\sqrt{q}. \eqno{\EQa1}$$
Our construction guarantees that $|g_i| \le \bigl(q^{n-i+1}+1\bigr)/2$
for $i=2$, \dots,~$n-1$, and it follows that 
$$|a_i| \le \frac{q}{2}\eqno{\EQai}$$
for $i=2$, \dots,~$n-1$ and that 
$|a_n| \le (q+1)/2$.
Also, note that if $a_{n-1} = q/2$ then
$$a_n = g_n - a_{n-1}(q+1) \le (q^2+1)/2 - (q/2)(q+1) < 0,$$
and likewise if $a_{n-1} = -q/2$ then $a_n > 0$.

If $(a_n,q)>1$ and $a_n<0$, replace $a_n$ with $a_n+(q+1)$ and replace
$a_{n-1}$ with $a_{n-1}-1$.  If $(a_n,q)>1$ and $a_n\ge0$, 
replace $a_n$ with $a_n-(q+1)$ and replace $a_{n-1}$ with $a_{n-1}+1$.
Note that whether or not we have changed the 
values of $a_n$ and $a_{n-1}$, equation~\EQnumpoints\ still
holds, equation~\EQa1\ still holds, equation~\EQai\ holds for $i=2$, \dots,~$n-2$, 
we have $(a_n,q) = 1$, and
we have $|a_{n-1}|\le(q+1)/2$ and $|a_n|\le q+1$.
The inequality for $a_{n-1}$ is the only non-obvious statement,
and it can be seen as follows:
in order for $|a_{n-1}|$ to be greater than $(q+1)/2$, 
either the original value for $a_{n-1}$ must have been $q/2$
and $a_n$ must have been positive and not coprime to $q$,
or the original value for $a_{n-1}$ must have been $-q/2$
and $a_n$ must have been negative and not coprime to $q$.
However, the comment at the end of the preceding paragraph shows 
that neither of these possibilities could have occurred.

We calculate that
$$\align
\left(\sum_{i=1}^{n-1} \left|\frac{a_i}{q^{i/2}}\right|\right) + \left|\frac{a_n}{2q^{n/2}}\right| 
         &\le B_q + \frac{q}{2q}  + \cdots + \frac{q}{2q^{(n-2)/2}}
                       +\frac{q+1}{2q^{(n-1)/2}} + \frac{q+1}{2q^{n/2}} \\
         &= \frac{1}{2}\frac{\sqrt{q}-2}{\sqrt{q}-1} 
               + \frac{1}{2}\left( 1 + \cdots + q^{-(n-1)/2} + q^{-n/2}\right)\\
         &< \frac{1}{2}\left(\frac{\sqrt{q}-2}{\sqrt{q}-1} + \frac{\sqrt{q}}{\sqrt{q}-1}\right)\\
         &= 1,
\endalign
$$
so by Lemma~\LLWeilpolys\ the polynomial
$$f=\bigl(x^{2n} + q^n\bigr) + a_1 \bigl(x^{2n-1}+q^{n-1}x\bigr) + \cdots 
    + a_{n-1}\bigl(x^{n+1}+qx^{n-1}\bigr) + a_n x^n $$
is an ordinary Weil $q$-polynomial, and by equation~\EQnumpoints\ we have $f(1) = m$.
This proves the theorem when $n>2$.

Now suppose $n=2$.  We begin as in the preceding case:
Pick $a_1$ so as to minimize the absolute value 
of $m - \bigl(q^2+1\bigr) - a_1\bigl(q+1\bigr)$, and then choose $a_2$ so that 
$m = \bigl(q^2+1\bigr) + a_1\bigl(q+1\bigr) + a_2$.
Note that the hypotheses of the theorem imply that
$$\Bigl| m - \bigl(q^2 + 1\bigr)\Bigr| < C_q\sqrt{q}(q+1)
= \Bigl(\bigl\lfloor B_q\sqrt{q}\bigr\rfloor + 1/2\Bigr) (q+1),$$
so again we have $|a_1|\le \bigl\lfloor B_q\sqrt{q}\bigr\rfloor \le B_q\sqrt{q}$.
Also, we see that $|a_2| \le (q+1)/2$.

Suppose $(a_2,q) = 1$.  Then Lemma~\LLWeilpolys, together with
the fact that 
$$\left|\frac{a_1}{\sqrt{q}}\right| + \left|\frac{a_2}{2q}\right|
   \le B_q + \frac{q+1}{4q} \
   = \frac{1}{2}\left(1 - \frac{1}{\sqrt{q}-1} + \frac{1}{2} + \frac{1}{2q}\right) < 1,$$
shows that $f = \bigl(x^4+q^2\bigr) + a_1\bigl(x^3+qx\bigr) + a_2 x^2$ is an ordinary Weil $q$-polynomial
with $f(1) = m$, and we are done.  

On the other hand, suppose that $(a_2,q)>1$.
Replace $a_2$ with $a_2+q+1$ and replace $a_1$ with $a_1-1$.  
We still have $m = (q^2+1) + a_1(q+1) + a_2$, 
and in addition we have $(q+1)/2\le a_2 \le 3(q+1)/2$
and $|a_1| < B_q\sqrt{q}+1 < \sqrt{q}$.  It follows that $a_1/\sqrt{q} \in [-1,1]$
and $a_2/2q\in[0,2]$.  By Example~\EEregion2, 
we have $(a_1/\sqrt{q},a_2/q)\in V_2$, so once again 
$f=\bigl(x^4+q^2\bigr) + a_1\bigl(x^3+qx\bigr) + a_2 x^2$ 
is an ordinary Weil $q$-polynomial with $f(1)=m$.
\myqed
\enddemo

\proclaim{Exercise \EXgrouporders}
Let $C_2 = 7\sqrt{2}/64$ and let $C_3 = 7\sqrt{3}/54$.
Suppose that $q=2$ or $q=3$ and that $n>1$ is an integer.
Show that if $m$ is an integer such that 
$\bigl|m - (q^n+1)\bigr| \le C_q q^{n-\frac{1}{2}}$
then there is an $n$-dimensional ordinary abelian variety $A$ over $\BF_q$
with $m = \#A(\BF_q)$.
\endproclaim

\demo{Hint}
Suppose $q=2$.  If $n\ge 7$, use the same argument as in the proof
of Theorem~\TTgrouporders, but start by taking
$a_1 = a_2 = a_3 = a_4 = 0$, and note that then $|a_5|\le 3$.
Check the cases $n<7$ by hand.
% checked 2/9/98
Similarly,
if $q=3$ and $n\ge 5$, use the same argument but start with $a_1 = a_2  = 0$,
and note that then $|a_3|\le 3$. 
Check the cases $n<5$ by hand.
% checked 2/9/98
\enddemo

\enddocument